\documentclass[seceqn,dvips]{arxbj}
\usepackage{dcolumn} 
\usepackage{cursive}
\usepackage{graphicx}

\aid{0}
\volume{15}
\issue{4}
\pubyear{2009}
\firstpage{1036}
\lastpage{1056}
\doi{10.3150/09-BEJ185}

\makeatletter
\newcolumntype{d}[1]{D{.}{.}{#1}}
\newcommand{\eqref}[1]{(\ref{#1})}
\newtheorem{lemma}{Lemma}
\newtheorem{theorem}{Theorem}
\newremark{C}{}
\makeatother

\begin{document}
\begin{frontmatter}

\title{Asymptotic optimal designs under long-range dependence error structure}
\runtitle{Designs under long-range dependence}

\begin{aug}
\author[a]{\fnms{Holger} \snm{Dette}\corref{}\thanksref{a}\ead[label=e1]{holger.dette@rub.de}},
\author[b]{\fnms{Nikolai} \snm{Leonenko}\thanksref{b}\ead[label=e2]{LeonenkoN@cf.ac.uk}},
\author[c]{\fnms{Andrey}~\snm{Pepelyshev}\thanksref{c}\ead[label=e3]{andrey@ap7236.spb.edu}}~\and\break
\author[d]{\fnms{Anatoly} \snm{Zhigljavsky}\thanksref{d}\ead[label=e4]{ZhigljavskyAA@cf.ac.uk}}
\runauthor{Dette, Leonenko, Pepelyshev and Zhigljavsky}
\address[a]{Ruhr-Universit\"{a}t Bochum, Fakult\"{a}t f\"{u}r
Mathematik, 44780 Bochum, Germany.\\\printead{e1}}
\address[b]{Cardiff University, School of Mathematics, Cardiff CF24
4AG, UK.\\\printead{e2}}
\address[c]{St. Petersburg State University, Department of Mathematics,
St. Petersburg, Russia.\\\printead{e3}}
\address[d]{Cardiff University, School of Mathematics, Cardiff CF24
4AG, UK.\\\printead{e4}}
\end{aug}

\received{\smonth{9} \syear{2008}}
\revised{\smonth{7} \syear{2009}}

%
\begin{abstract}
We discuss the optimal design problem in regression models with
long-range dependence error structure. Asymptotic optimal designs are derived
and it is demonstrated that these designs depend only indirectly on the
correlation function.
Several examples are investigated to illustrate the theory.
Finally, the optimal designs are compared with
asymptotic optimal designs which were derived by Bickel and Herzberg
[\textit{Ann. Statist.} \textbf{7} (1979) 77--95]
for regression models with short-range dependent error.
\end{abstract}

%

\end{frontmatter}

\section{Introduction}

Consider the common linear regression model
%
\begin{eqnarray}\label{1.1}
y(t)=\theta_1f_1(t)+\cdots+\theta_pf_p(t)+\varepsilon(t),
\end{eqnarray}
where $f_1(t), \ldots, f_p(t)$ are known functions,
$\varepsilon(t)$ is a random error, $\theta_1, \ldots, \theta_p$
denote the unknown parameters and $t$ is the explanatory variable.
We assume that $N$ observations, say $y_1,\ldots, y_N$, can be
taken at experimental conditions $-T\le t_1\le\cdots\le t_N\le T$ to estimate the parameters in the
linear regression model (\ref{1.1}). If an appropriate estimate,
say $ \hat\theta= ( \hat\theta_1, \ldots, \hat\theta_p)^T$, has
been chosen, an optimal design minimizes a function of the variance-covariance
matrix of this estimate, which is called an optimality criterion
(see, for example, Silvey (\citeyear{Silvey1980}) or Pukelsheim
(\citeyear{Pukelsheim1993})).

Under the assumption of uncorrelated observations, optimal designs have
been studied by numerous authors (see the two books cited above and the
textbooks of Fedorov (\citeyear{Fedorov1972}),
P\'{a}zman (\citeyear{Pazman1986}) and
Atkinson and Donev (\citeyear{AtkinsonDonev1992})).
However, fewer results are available for dependent
observations, although this problem is of particular interest because
in many applications, the variable $t$ in the regression model
(\ref{1.1}) represents time and all observations correspond to one
subject. The reason for this is that optimal experimental designs for
regression models with correlated observations have a very complicated
structure and are difficult to find, even in simple cases. Because
explicit solutions are rarely available, an asymptotic theory was
developed by Sacks and Ylvisaker (\citeyear{SacksYlvisaker1966}, \citeyear{SacksYlvisaker1968}).
In the Sacks--Ylvisaker approach, the design set is fixed and the number of design points in
this set tends to infinity. As a result of this assumption, the design
points become too close to each other and the corresponding asymptotic
optimal designs depend only on the behavior of the correlation function
in a neighborhood of the point 0.

Bickel and Herzberg (\citeyear{BickelHerzberg1979}) and
Bickel, Herzberg and Schilling (\citeyear{BickelHerzbergSchilling1981})
considered a~different model,
where the design interval expands proportionally to the number of
observation points and
the correlation structure of errors is not used for the construction of
the least-squares estimate.
The variance-covariance matrix of the estimate $\hat\theta$ is of
order~$\mathrm{O}(1)$ in the model considered by Sacks and Ylvisaker (\citeyear{SacksYlvisaker1966})
and of order $1/N$ in the model discussed by Bickel and Herzberg (\citeyear{BickelHerzberg1979}).
Therefore, the approach of Bickel and Herzberg makes the optimal
designs derived for the dependent and independent cases more comparable.
These authors assumed that the observations
in model (\ref{1.1}) have a correlation structure corresponding
to a non-degenerate stationary process with
short-range dependence, where
a correlation function
$\rho$ satisfies $\rho(t) =\mathrm{o}(1/t)$ if $t \to\infty$.
As examples, in Bickel and Herzberg (\citeyear{BickelHerzberg1979}) and
Bickel, Herzberg and Schilling (\citeyear{BickelHerzbergSchilling1981}),
asymptotic
optimal designs are derived for the linear regression model with and without
intercept and for the location model.

The purpose of the present paper is to extend the Bickel--Herzberg
approach to the case
of a stronger dependence
of the errors in the linear regression model (\ref{1.1}), which
corresponds to an
error process with long-range dependence. Long-range dependence is
observed in many applications, including hydrology,
geophysics, turbulence,
diffusion,
economics and finance.
The phenomenon was already observed by Pearson (\citeyear{Pearson1902}) in
astronomy and by Smith (\citeyear{Smith1938})
in agriculture. Further examples where long-range dependence is
discussed can be found in Granger (\citeyear{Granger1980}),
Mandelbrot (\citeyear{Mandelbrot1973}),
Porter-Hudak (\citeyear{Porterhudak1990}),
Beran, Sherman, Taqqu and Willinger (\citeyear{BeranEtAl1992}),
Barndorff-Nielsen \textit{et al.}~(\citeyear{BarndnielsenJensenSorensen1990}),
Beran (\citeyear{Beran1992}), Metzler \textit{et al.}~(\citeyear{MetzlerBarkaiKlafter1999}),
among many others. The interested reader is referred to the books of
Beran (\citeyear{Beran1994}) and
Doukhan \textit{et al.}~(\citeyear{DoukhanOppenheimTaqqu2003}), which contain a good
description of the basic properties
of long-range dependence processes and
an extensive bibliography on this subject.

Most of the literature considers the estimation problem but -- to the
best knowledge of the authors -- design problems for regression models with
long-range dependence error structure have not been considered thus
far. In Section~\ref{sec2}, we introduce the basic terminology and describe the optimal
design problem. Our main results are given in Section~\ref{sec3}, where we
derive an asymptotic expression for the variance-covariance matrix,
the basis for the construction of optimal designs in the regression
model (\ref{1.1}) with a long-range dependent error structure.
These results are different from the findings of
K\"{u}nsch, Beran and Hampel (\citeyear{KunschBeranHampel1993}), who considered random explanatory variables.
Finally, in Section~\ref{sec4} several asymptotic optimal designs
are derived for the linear regression model and compared with the results
obtained by Bickel and Herzberg (\citeyear{BickelHerzberg1979})
 under the assumption of a short-range
error structure.

\section{Optimal designs for dependent observations}\label{sec2}

Consider the linear regression model (\ref{1.1}), where the error
process $\varepsilon(t)$
is the second-order process with
%
\begin{eqnarray}\label{2.1}
\mathbf{E}\varepsilon(t)=0,\qquad
\mathbf{E}\varepsilon(t)\varepsilon(s)=\sigma^2 r(t,s),
\end{eqnarray}
and assume that
\renewcommand{\theC}{\textup{(C1)}}
\begin{C}\label{C1}
The regression
functions $f_1(t),\ldots,f_p(t)$ are linearly independent and bounded on
the interval $[-T,T]$ and satisfy
a first order Lipschitz condition, that is,
$|f_i(t)-f_i(s)|\le M|t-s|$ and $|f_i(t)|\le M$
for all $t,s\in[-T,T]$, $i=1,\ldots,p$.
\end{C}

Following Bickel and Herzberg (\citeyear{BickelHerzberg1979}), we assume that
$\varepsilon(t)=\varepsilon^{(1)}(t)+\varepsilon^{(2)}(t)$, where~$\varepsilon^{(1)}(t)$
denotes a stationary process with correlation function $\rho(t)$ and
$\varepsilon^{(2)}(t)$ is white noise. Consequently, we obtain
%
\begin{eqnarray}\label{eq:corr_fun_r}
r(t,s)=\gamma\rho(t-s)+(1-\gamma)\delta_{t,s},
\end{eqnarray}
where $\delta$ is the Kronecker symbol. If $N$ observations, say
$y=(y_1, \ldots,y_N)^T$, are
available and the form of the correlation function is known, then the
vector of parameters can
be estimated by the weighted least squares, that is,
$\hat\theta= (X^T\Sigma^{-1}X)^{-1}X^T\Sigma^{-1}y$
with $X^T=(f_i(t_j))_{i=1,\ldots,p}^{j=1,\ldots,N}$,
and the variance-covariance matrix of this estimate is given by
%
\begin{eqnarray}\label{covwlse}
\mathbf{D} (\hat\theta) =\sigma^2(X^T\Sigma^{-1}X)^{-1}
\end{eqnarray}
with
$\Sigma=(\gamma\rho(t_i-t_j)+(1- \gamma)\delta_{i,j})_{i,j}$,
$i,j=1,\ldots,N$.
However, in most applications, knowledge about the correlation
structure is not available
and the unweighted least-squares estimate
$ \tilde\theta= (X^TX)^{-1}X^Ty$ is used. For this estimate, the
variance-covariance matrix is given by
%
\begin{eqnarray}
\mathbf{D} (\tilde\theta) =\sigma^2(X^TX)^{-1} X^T\Sigma X (X^TX)^{-1}.
\label{eq:lse}
\end{eqnarray}
An experimental design $\xi=\{t_1,\ldots,t_N\}$ is a vector of
$N$ points in the interval $[-T,T]$, which defines the time points
or experimental conditions where observations are taken. Optimal
designs minimize a functional of the variance-covariance matrix of
the weighted or unweighted least-squares estimate. Following
Bickel and Herzberg (\citeyear{BickelHerzberg1979}), we consider a~correlation function which
depends on the sample size $N$ and is of the form
\mbox{$\rho_N(t)=\rho(Nt)$}, where the function $\rho$ satisfies $\rho(t)
\to0 $ if $t\to\infty$; this corresponds to expanding the\vadjust{\goodbreak}
interval as the number of observations grows. The standard
least-squares estimate is considered in the following discussion
because when computing this estimate, the
form of the correlation function $\rho(t)$ is not used. Despite
this, the least-squares estimate often has good properties compared
to the best linear unbiased estimate; see, for example,
Adenstedt (\citeyear{Adenstedt1974}),
Samarow and Taqqu (\citeyear{SamarovTaqqu1988}), Yajima (\citeyear{Yajima1988}, \citeyear{Yajima1991}) and
Beran (\citeyear{Beran1994}), page~179, among others.
For some results regarding
nonlinear regression, the reader is referred to
Ivanov and Leonenko (\citeyear{IvanovLeonenko2004}, \citeyear{IvanovLeonenko2008}).

For our asymptotic investigations, we consider a sequence of
designs $\xi_N=  \{t_{1N}, \ldots,  t_{NN}\}$
which is generated using a continuous non-decreasing function
%
\begin{eqnarray}\label{dens}
a{}\dvtx{}[0,1]\to[-T,T]
\end{eqnarray}
by
%
\begin{eqnarray}\label{design}
t_{iN}=a\bigl((i-1)/(N-1)\bigr),\qquad i=1,\ldots,N,
\end{eqnarray}
where the function $a(u)$ is the
inverse of a distribution function. Note that the function~$a$ is
obtained as the
weak limit of $\xi_N$ as $N\to\infty$. The
equally spaced design corresponds to the choice $a(u)=(2u-1)T$ ($u \in
[0,1]$); changing the function $a$ yields different types of the
design. For example,
the choice $a^{-1}(x)=(x^3+1)/2$ yields designs which are more
concentrated in the interval
$[-1/2, 1/2]$.
We assume several regularity conditions on the function $a$, which are
required for the asymptotic results which are to follow.
More precisely:\looseness=-1

\renewcommand{\theC}{\textup{(C2)}}
\begin{C}\label{C2}
Let $a(u)$ be twice differentiable and assume that
there exists a positive constant $M<\infty$ such that for all $ u\in(0,1)$,
%
\begin{eqnarray}\label{bound}
\frac{1}{M}\le a'(u) \le M,\qquad |a''(u)|\le M.
\end{eqnarray}
\end{C}

\renewcommand{\theC}{\textup{(C3)}}
\begin{C}\label{C3}
The correlation function $\rho$ is differentiable with bounded
derivative, that is, $|\rho'(t)|\le M$, $t\in(0,\infty)$ and
$\rho'(t)\le0$ for sufficiently large $t$.
\end{C}

The last assumption implies that $\rho(t)$ is nonnegative for sufficiently large $t$. In contrast to Bickel and Herzberg (\citeyear{BickelHerzberg1979})
we assume that
%
\begin{eqnarray}\label{inte}
\int_0^\infty|\rho(t)|\,\mathrm{d}t=\infty
\end{eqnarray}
and this assumption corresponds to
the long-range dependence of the observations. Note that in this case
it follows that
\begin{eqnarray*}
\int_0^\infty|\rho(t)|\,\mathrm{d}t=\sum_{k=0}^\infty|\rho(k)|=\infty,
\end{eqnarray*}
where $\rho(k)=\operatorname{cov}(\varepsilon^{(1)}(t), \varepsilon^{(1)}(t+k))$.
The correlation function of a stationary process with long-range dependence
can be written as
%
\begin{eqnarray}\label{eq:rho=L/ta}
\rho_\alpha(t)=\frac{L(t)}{|t|^\alpha},\qquad |t|\to\infty,\vadjust{\goodbreak}
\end{eqnarray}
where $0<\alpha\le1$ and $L(t)$ is a slowly varying function (SVF)
for large $t$
[Doukhan \textit{et al.} (\citeyear{DoukhanOppenheimTaqqu2003})],
and satisfies
\begin{eqnarray*}
\rho_\alpha(t)=\mathrm{O}(1/|t|^\alpha),\qquad |t|\to\infty.
\end{eqnarray*}
In this case, we will say that $\rho_\alpha(t)$ belongs to
the SVF family.

\section{Main results}\label{sec3}

First, we introduce two parametric families of correlation functions
which are important in applications.

\begin{figure}

\includegraphics{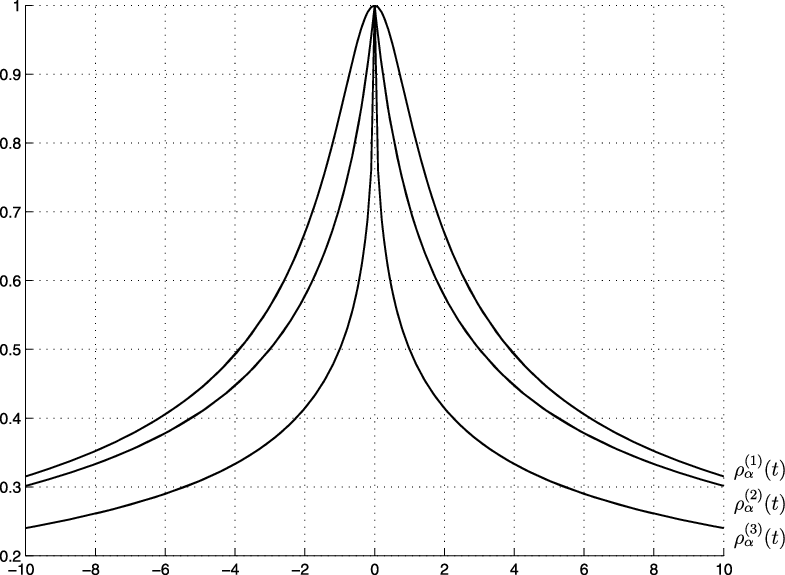}

\caption{The three correlation functions, where $\alpha=
0.5$.}\label{cauchy}
\end{figure}

The correlation function $\rho_\alpha(t)$
belongs to the
Cauchy family if it is defined by
%
\begin{eqnarray}\label{eq:cauchy_fam}
\rho_\alpha(t)=\frac{1}{(1+|t|^\beta)^{\alpha/\beta}},
\end{eqnarray}
where $\beta>0$, $0<\alpha\le1$ [see Gneiting (\citeyear{Gneiting2000}),
Anh \textit{et al.}~(\citeyear{AnhKnopovaLeonenko2004}),
Barndorff-Nielsen and Leonenko (\citeyear{BarndnielsenLeonenko2005})]. This family includes
\begin{eqnarray*}
\rho^{(1)}_\alpha(t)=\frac{1}{(1+|t|^2)^{\alpha/2}},\qquad
\rho^{(2)}_\alpha(t)=\frac{1}{1+|t|^\alpha},\qquad
\rho^{(3)}_\alpha(t)=\frac{1}{(1+|t|)^\alpha},
\end{eqnarray*}
which have a totally different shape in a neighbourhood of the point $t=0$
but the same asymptotic behavior for large $t$ (see Figure~\ref{cauchy}).
These three functions are known as the characteristic functions of the
symmetric Bessel distribution,
Linnik distribution and symmetric generalized Linnik distribution, respectively.

The correlation function $\rho_\alpha(t)$
belongs to the Mittag-Leffler family if it is defined\break by
%
\begin{eqnarray}\label{eq:mittag-leffler_fam}
\rho_\alpha(t)=E_{\nu,\beta}(-|t|^\alpha),\qquad
E_{\nu,\beta}(-t)=\Gamma(\beta)\sum_{k=0}^\infty
\frac{(-t)^k}{\Gamma(\nu k+\beta)}, \qquad t > 0,
\end{eqnarray}
where $0<\alpha\le1$, $0<\nu\le1$, $\beta\ge\nu$ (see Schneider (\citeyear{Schneider1996}),
Barndorff-Nielsen and Leonenko (\citeyear{BarndnielsenLeonenko2005})).
This family is a smooth interpolation of long-range dependence
($ 0<\alpha\leq1$, $\beta=1, 0 < \nu< 1$)
and short-range dependence $(\nu=1,0< \alpha\leq1, \beta=1)$.
Note that the case $\nu=1,\beta=1,\alpha=1$ corresponds to ordinary
diffusion $\varrho_\alpha(t)= \mathrm{e}^{-|t|}$,
which is the correlation function of a Markovian Ornstein--Uhlenbeck
process. On the other hand,
the case $0<\nu<1,\beta=1,0<\alpha<1$ corresponds to subdiffusion or
slow diffusion
(see Metzler and Klafter (\citeyear{MetzlerKlafter2000})).
In particular,
\begin{eqnarray*}
E_{1,1}(-t)&=&\mathrm{e}^{-t},\qquad
E_{1,2}(-t)=(1-\mathrm{e}^{-t})/t,\qquad
E_{1,3}(-t)=2(\mathrm{e}^{-t}-1+t)/t^2,
\\[3pt]
E_{1/2,1}(-t)&=&\mathrm{e}^{t^2}\biggl(1-\frac{2}{\sqrt{\curpi}}
\int_0^t\mathrm{e}^{-u^2}\,\mathrm{d}u\biggr).
\end{eqnarray*}
In the following discussion, we derive optimal designs for the three
families of correlation functions, which are given by
\eqref{eq:rho=L/ta}, \eqref{eq:cauchy_fam} and \eqref{eq:mittag-leffler_fam}.
The function $Q(t)=\sum_{j=1}^\infty\rho(j t)$\vspace*{2pt}
plays an important role in the asymptotic analysis by Bickel and Herzberg (\citeyear{BickelHerzberg1979}),
but in the case of long-range dependence, this function is infinite.
For an asymptotic analysis under long-range dependence, we
introduce the function
%
\begin{eqnarray}
\label{qalph}
Q_\alpha(t)=\lim_{N\to\infty}\frac{1}{d_\alpha(N)}\sum_{j=1}^N\rho_\alpha(j t),
\end{eqnarray}
where the normalizing sequence is given by
\begin{eqnarray*}
d_\alpha(N)=
\cases{
N^{1-\alpha},
&\quad\mbox{if }$\alpha<1$\mbox{ and }$\rho_\alpha$\mbox{ has the form }\eqref{eq:cauchy_fam}\mbox{ or }\eqref{eq:mittag-leffler_fam},
\cr
\ln N, &\quad\mbox{if }$\alpha=1$\mbox{ and }$\rho_\alpha$\mbox{ has the form }\eqref{eq:cauchy_fam}\mbox{ or }\eqref{eq:mittag-leffler_fam},
\cr
L(N) N^{1-\alpha}, &\quad\mbox{if }$\alpha<1$\mbox{ and }$\rho_\alpha$\mbox{ has the form }\eqref{eq:rho=L/ta},
\cr
L(N) \ln N, &\quad\mbox{if }$\alpha=1$\mbox{ and }$\rho_\alpha$\mbox{ has the form }\eqref{eq:rho=L/ta},
}
\end{eqnarray*}
and show in Lemma~\ref{lemma1} below that the function $Q_\alpha(t)$ is well
defined.



\begin{lemma}\label{lemma1}
If the correlation function $\rho_\alpha(t)$ belongs to the Cauchy,
SVF family or to the Mittag-Leffler family with
$0< \alpha\leq1, 0 < \nu\leq1, \nu\leq\beta, (\nu,\beta) \neq
(1,1)$, then the limit in~(\ref{qalph}) exists and is given by
\begin{eqnarray*}
Q_\alpha(t)=
\cases{
\displaystyle\frac{c}{(1-\alpha)|t|^\alpha},&\quad $0<\alpha<1$,\cr
\displaystyle\frac{c}{|t|},&\quad$\alpha=1$,
}
\end{eqnarray*}
where
\begin{eqnarray*}
c=
\cases{
\displaystyle\frac{\Gamma(\beta)}{\Gamma(\beta-\nu)},
&\quad\mbox{\textup{if} } $\rho_\alpha(t)$\mbox{ \textup{belongs to the Mittag-Leffler family}},
\cr
1,&\quad\mbox{\textup{otherwise}}.
}
\end{eqnarray*}
\end{lemma}

\begin{pf}
Define the function
\begin{eqnarray*}
Q_{\alpha,N}(t)=\frac{1}{d_\alpha(N)} \sum_{j=1}^N\rho_\alpha(j t)
\end{eqnarray*}
and assume that the correlation function $\rho_\alpha(t)$
is an element of the Cauchy family.
Since the function $\rho_\alpha(t)$ defined in \eqref{eq:cauchy_fam}
is positive and decreasing for $0<\alpha<1$, we have
\begin{eqnarray*}
Q_{\alpha,N}(t)&=&\frac{1}{N^{1-\alpha}}\int_0^N
\frac{1}{(1+|st|^\beta)^{\alpha/\beta}}\, \mathrm{d}s +\mathrm{O}\biggl(\frac{1}{N^{1-\alpha}}\biggr)
\\
&=&\frac{1}{N^{1-\alpha}}N\int_0^N\frac{\mathrm{d}(s/N)}{N^\alpha(1/N^\beta+|st/N|^\beta)^{\alpha/\beta}}
+\mathrm{O}\biggl(\frac{1}{N^{1-\alpha}}\biggr)
\\
&=&\int_0^1\frac{\mathrm{d}v}{(|vt|^\beta)^{\alpha/\beta}}+\mathrm{O}\biggl(\frac{1}{N^{\alpha-\alpha^2}}\biggr)
=
\frac{1}{(1-\alpha)|t|^\alpha}+\mathrm{O}\biggl(\frac{1}{N^{\alpha-\alpha^2}}\biggr)
\\
&=& Q_\alpha(t) +\mathrm{O}\biggl(\frac{1}{N^{\alpha-\alpha^2}}\biggr).
\end{eqnarray*}
For $\alpha=1$, we obtain
\begin{eqnarray*}
Q_1(t)&=&\lim_N\frac{1}{\ln N}\sum_{j=1}^N\frac{1}{(1+|j t|^\beta)^{1/\beta}}
=
\lim_N\frac{1}{|t|\ln N}\int_0^N \frac{1}{(1+|st|^\beta)^{1/\beta}}\, \mathrm{d} (st)\\
&=& \lim_N\frac{1}{|t|\ln N}\int_0^{Nt} \frac{1}{1+|v|}\, \mathrm{d} v
=\frac{1}{|t|},
\end{eqnarray*}
which completes the proof for the case where $\rho_\alpha(t)$ belongs
to the Cauchy family.

Now assume that the correlation function $\rho_\alpha(t)$
is an element of the Mittag-Leffler family. For the sake of brevity, we
only consider the case
$\beta=1, 0 < \alpha\leq1, 0 <\nu<1$, all other cases being treated
similarly.
Since
\[
E_{\nu,1}(-|t|^\alpha)\sim\frac{1}{|t|^\alpha\Gamma(1-\nu)}
\]
as $t\to\infty$ (see, for example, formula (3.17) in
Schneider (\citeyear{Schneider1996})), we have, for $0 < \alpha<1$,
\begin{eqnarray*}
Q_\alpha(t)=\lim_{N \to\infty}\frac{1}{N^{1-\alpha}}\sum_{j=1}^N
E_{\nu,1}(-|jt|^\alpha)
=\frac{1}{(1-\alpha)\Gamma(1-\nu)|t|^\alpha}.
\end{eqnarray*}
Observing
that
\[
E_{\nu,\beta}(-|t|)\sim\frac{\Gamma(\beta)}{|t|\Gamma(\beta-\nu)}
\]
for $t\to\infty$ (see Djrbashian (\citeyear{Djrbashian1993})), we obtain, for $\alpha=1$,
\begin{eqnarray*}
Q_1(t)=\lim_{N \to\infty}\frac{1}{\ln N}\sum_{j=1}^N E_{\nu,\beta}(-|jt|)
=\frac{\Gamma(\beta)}{|t|\Gamma(\beta-\nu)}.
\end{eqnarray*}
Finally, assume that the correlation function $\rho_\alpha(t)$ is an
element of the SVF family.
We then obtain
\begin{eqnarray*}
Q_{\alpha}(t)&=&\lim_{N \to\infty}\frac{1}{L(N)N^{1-\alpha}}
\int_0^N \frac{L(st)}{|st|^{\alpha}}\, \mathrm{d}s
\\
&=&\lim_{N \to\infty}\frac{1}{L(N)N^{1-\alpha}}N\int_0^N
\frac{L(N t s/N)\,\mathrm{d}(s/N)}{N^\alpha|st/N|^{\alpha}}
\\
&=&\lim_{N \to\infty}\int_0^1\frac{L(N t v)\,\mathrm{d}v}{L(N)|vt|^{\alpha}}
=\int_0^1\frac{\mathrm{d}v}{|vt|^{\alpha}}=
\frac{1}{(1-\alpha)|t|^\alpha},
\end{eqnarray*}
where we have used Theorem 2.6 from Seneta (\citeyear{Seneta1976}) in the last line.
For $\alpha=1$, we have
\begin{eqnarray*}
Q_1(t)&=&\lim_{N \to\infty}\frac{1}{L(N)\ln N}\int_{1}^N\frac{L(st)}{|st|}\,\mathrm{d}s =
\lim_N\frac{1}{\ln N}\int_{1}^N\frac{L(st)/L(N)}{st}\,\mathrm{d}s \\
&=&\lim_{N \to\infty}\frac{1}{\ln N}\int_{1}^N\frac{1}{|st|}\,\mathrm{d}s+
\lim_N\frac{1}{\ln N}\int_{1}^N\frac{L(st)/L(N)-1}{|st|}\,\mathrm{d}s\\
&=&\frac{1}{|t|}+\lim_N\frac{1}{\ln N}\int_{1/N}^1\frac{L(Nvt)/L(N)-1}{|vt|}\,\mathrm{d}v=\frac{1}{|t|},
\end{eqnarray*}
which completes proof of Lemma~\ref{lemma1}.
\end{pf}

Next, we find a comfortable asymptotic representation for the main term
in the
variance-covariance matrix of the least-squares estimates.

\begin{lemma}\label{lemma2}
Assume that the correlation function $\rho_\alpha(t)$ belongs to the Cauchy,
Mittag-Leffler or SVF family, such that
%
\begin{equation}\label{ass}
\int_0^1 Q_\alpha(a'(t))\,\mathrm{d}t<\infty,
\end{equation}
and that the regularity conditions \textup{\ref{C1}}--\textup{\ref{C3}} in
Sections~\ref{sec2}
and~\ref{sec3} are satisfied. We have
\begin{eqnarray*}
\frac{1}{d_\alpha(N) N}
\sum_{i\neq j}f_s(t_{iN})f_r(t_{jN})\rho_\alpha\bigl(N(t_{jN}-t_{iN})\bigr)
=2\int_0^1 f_s(a(u))f_r(a(u))Q_\alpha(a'(u))\,\mathrm{d}u+\mathrm{o}(1)
\end{eqnarray*}
as $N\to\infty$ for all $s,r=1, \ldots,p$, $0<\alpha\le1$.
\end{lemma}


\begin{pf}
We only give a proof for the correlation function from the Cauchy
family and $0<\alpha<1$,
the proof for the other cases being similar. We use the notation
$f=f_s$, $g=f_r$, $\rho= \rho_\alpha$ and the decomposition
\begin{eqnarray*}
N^{\alpha-2}\sum_{i\neq j}f(t_{iN})g(t_{jN})\rho\bigl(N(t_{jN}-t_{iN})\bigr)=S_1+S_2,
\end{eqnarray*}
where
%
\begin{eqnarray}\label{s1}
S_1&=&2N^{\alpha-2}\sum_{i=1}^N f(t_{iN})g(t_{iN})\sum_{j=i+1}^N\rho
\bigl(N(t_{jN}-t_{iN})\bigr),
\\\label{s2}
S_2 &=&2N^{\alpha-2}\sum_{i=1}^Nf(t_{iN})\sum_{j=i+1}^N\bigl(g(t_{jN})-g(t_{iN})\bigr)\rho\bigl(N(t_{jN}-t_{iN})\bigr).
\end{eqnarray}
With the notation $i_N=(i-1)/(N-1)$, we obtain from the
differentiability of the functions
$a$ and $\rho$
\begin{eqnarray*}
\rho\bigl(N(t_{jN}-t_{iN})\bigr)
=\rho\bigl(N\bigl(a(j_N)-a(i_N)\bigr)\bigr)
=\rho\bigl(a'(i_N)(j-i)\bigr)+\nu\frac{(j-i)^2}{N-1},
\end{eqnarray*}
where $|\nu|\le M^2/2$.
Let $r_N$ denote a sequence such that $r_N\to\infty$ slowly as
$\mathrm{o}(N^{(1-\alpha)/3})$ and
consider the cases $i \le r_N$ and $i > r_N$ in (\ref{s1}) and (\ref{s2}) separately. Note that
\begin{eqnarray*}
\Bigg|\sum_{j=i+r_N}^N \rho\bigl(N(t_{jN}-t_{iN})\bigr)\Bigg|
&=& \sum_{j=i+r_N}^N \rho\bigl(N\bigl(a(j_N)-a(i_N)\bigr)\bigr)
\\
&\le&\tilde M\sum_{j=i+r_N}^N \rho\bigl((j-i)/M\bigr) \le\tilde M\sum_{k=r_N}^\infty\rho(k/M)
=\mathrm{o}(N^{1-\alpha})
\end{eqnarray*}
as $N\to\infty$ uniformly with respect to $j$,
where $\tilde M$ is a constant and where we have used the fact that the function
$a'(u)$ is bounded from below and Lemma~\ref{lemma1}. Similarly,
we obtain
\begin{eqnarray*}
&&\Bigg|\sum_{j=i+r_N}^N \bigl(g(t_{jN})-g(t_{iN})\bigr)\rho\bigl(N(t_{jN}-t_{iN})\bigr)\Bigg|
\\
&&\quad\le2MT\sum_{j=i+r_N}^\infty
\big|\rho\bigl(N(t_{jN}-t_{iN})\bigr)\big|
=\mathrm{o}(N^{1-\alpha})
\end{eqnarray*}
as $N\to\infty$ uniformly with respect to $j$
because the function $g$ is bounded. This implies
that
%
\begin{eqnarray}\label{s1a}
S_1&=&2N^{\alpha-2}\sum_{i=1}^N f(t_{iN})g(t_{iN})\sum_{j=i+1}^{i+r_N}\rho\bigl(N(t_{jN}-t_{iN})\bigr) + \mathrm{o}(1),
\\\label{s2a}
S_2 &=&2N^{\alpha-2}\sum_{i=1}^Nf(t_{iN})\sum_{j=i+1}^{i+r_N}
\bigl(g(t_{jN})-g(t_{iN})\bigr)\rho\bigl(N(t_{jN}-t_{iN})\bigr)
+ \mathrm{o}(1)
\end{eqnarray}
as $N\to\infty$.
For the first term on the right-hand side of (\ref{s2a}), we
obtain the estimate
\begin{eqnarray*}
\tilde S_2 &=& N^{\alpha-2}
\Bigg|\sum_{i=1}^N f(t_{iN})
\sum_{j=i+1}^{i+r_N}\bigl(g(t_{jN})-g(t_{iN})\bigr)\rho\bigl(N(t_{jN}-t_{iN})\bigr)\Bigg|
\\
& \le&
2N^{\alpha-1}M^2T\sum_{j=i+1}^{i+r_N}
\big|\rho\bigl(N(t_{jN}-t_{iN})\bigr)\big|
\\
&\le&
2N^{\alpha-1}M^2T\sum_{j=i+1}^{i+r_N}
\biggl(\big|\rho\bigl(a'(i_N)(j-i)\bigr)\big|+M^2\frac{(j-i)^2}{N-1}\biggr)
\\
& \le& 2N^{\alpha-1} M^2T(Mr_N +M^2r_N^3/N)=\mathrm{o}(1)
\end{eqnarray*}
as $N\to\infty$, while the dominating term on the right-hand side of
(\ref{s1a}) is given
by
\begin{eqnarray*}
\tilde S_1 &=& N^{\alpha-2}\sum_{i=1}^N f(t_{iN})g(t_{iN})
\sum_{j=i+1}^{i+r_N}\rho\bigl(N(t_{jN}-t_{iN})\bigr)
\\
& =&
N^{\alpha-2}\sum_{i=1}^N f(t_{iN})g(t_{iN})\sum_{j=i+1}^{i+r_N}
\rho\bigl(a'(i_N)(j-i)\bigr)+\mathrm{o}(1)
\\
&=& N^{-1}\sum_{i=1}^N f(t_{iN})g(t_{iN})Q_\alpha(a'(i_N))+\mathrm{o}(1)
= \int_0^1 f(a(u))g(a(u))Q_\alpha(a'(u))\,\mathrm{d}u+\mathrm{o}(1)
\end{eqnarray*}
as $N\to\infty$, which proves the assertion of
Lemma~\ref{lemma2}.\vadjust{\goodbreak}
\end{pf}

\begin{theorem}
Let the correlation function $\rho_\alpha(t)$ be an element of the Cauchy,
Mittag-Leffler or SVF family. If (\ref{ass})
and the regularity assumptions \textup{\ref{C1}}--\textup{\ref{C3}}
stated in Sections~\ref{sec2} and~\ref{sec3} are satisfied, then we obtain for the
variance-covariance matrix of the
least-squares estimate defined in (\ref{eq:lse})%
\begin{eqnarray*}
\sigma^2\frac{N}{d_\alpha(N)}\mathbf{D} (\tilde\theta) =
2\gamma W^{-1}(a)R_\alpha(a)W^{-1}(a)+\mathrm{O}\bigl(1/d_\alpha(N)\bigr),
\end{eqnarray*}
where the matrices $W$ and $R_\alpha$ are given by
\begin{eqnarray*}
W(a)&=&\biggl(\int_0^1 f_i(a(u)) f_j(a(u))\, \mathrm{d}u\biggr)_{i,j=1}^p,
\\
R_\alpha(a)&=&\biggl(\int_0^1 f_i(a(u)) f_j(a(u))
Q_\alpha(a'(u))\,\mathrm{d}u\biggr)_{i,j=1}^p.]
\end{eqnarray*}
\end{theorem}

\begin{pf}
In view of \eqref{eq:corr_fun_r},
we obtain that
\begin{eqnarray*}
X^T\Sigma X
= \Biggl( \gamma\sum_{i\neq j}f_k(t_{iN})f_l(t_{jN})\rho_\alpha\bigl(N(t_{jN}-t_{iN})\bigr)+
\sum_{i=1}^N f_k(t_{iN})f_l(t_{jN})\Biggr)_{k,l=1}^p,
\end{eqnarray*}
where $X^T=(f_i(t_{jN}))_{i=1,\ldots,p}^{j=1,\ldots,N}$ and
$t_{iN}=a((i-1)/(N-1)), i=1,\ldots,N$.
An application of Lemma~\ref{lemma2} yields
\begin{eqnarray*}
\frac{X^TX}{N}= W(a)+\mathrm{O}\biggl(\frac{1}{N}\biggr),\qquad
\frac{X^T\Sigma X}{d_\alpha(N)N}
= 2\gamma R_\alpha(a)+ \mathrm{O}\biggl(\frac{1}{d_\alpha(N)}\biggr).
\end{eqnarray*}
The assertion of the theorem now follows by inserting these limits into
\eqref{eq:lse}.
\end{pf}

Note that the constant $ \gamma$ only
appears as a factor in the asymptotic variance-covariance matrices of
the least-squares estimate. Because most optimality
criteria are positively homogeneous (see, for example,
Pukelsheim (\citeyear{Pukelsheim1993})), it is reasonable to consider
the matrix
\begin{eqnarray*}
W^{-1}(a)R_\alpha(a)W^{-1}(a),
\end{eqnarray*}
which is proportional to the asymptotic variance-covariance matrix of
the least-squares estimate.
Moreover, if the function $a$ corresponds to a continuous distribution
with a~density, say $\phi$, then $a'(t)=1/\phi(t)$ and the
asymptotic variance-covariance matrix of the least-squares estimate is
proportional to the matrix
\begin{eqnarray*}
\Psi_\alpha(\phi)=W^{-1}(\phi)R_\alpha(\phi)W^{-1}(\phi),
\end{eqnarray*}
where the matrices $W(\phi)$ and $R_\alpha(\phi)$ are given by
\begin{eqnarray*}
W(\phi) &=&\biggl(\int_{-T}^T f_i(t) f_j(t)\phi(t)\, \mathrm{d}t\biggr)_{i,j=1,\ldots,p},
\\R_\alpha(\phi)
&=&\biggl(\int_{-T}^T f_i(t) f_j(t) Q_\alpha\bigl(1/\phi(t)\bigr)\phi(t)\,\mathrm{d}t\biggr)_{i,j=1,\ldots,p}
\\
&
=&\frac{c}{1-\alpha}\biggl(\int_{-T}^T f_i(t) f_j(t) \phi^{1+\alpha}(t)\,\mathrm{d}t\biggr)_{i,j=1,\ldots,p}
\end{eqnarray*}
and we have used the representation $Q_\alpha(t)=c/((1-\alpha)|t|^\alpha)$ for the last identity.
An (asymptotic) optimal design minimizes an appropriate function of the matrix
$\Psi_\alpha(\phi)$ (for classical least-squares
estimation).
Note that under long-range dependence, the variance-covariance matrix
of the least-squares estimate converges slower to zero
than in the case of independent or short-range dependent errors.
In the case of short-range dependence, no other normalization is
necessary, apart
from normalizing the variance-covariance matrix. Under long-range
dependence, an additional factor $d_{\alpha}(N)/N$ is
needed. Moreover, it is worthwhile to
note that under long-range dependence, the asymptotic
variance-covariance matrix is fully determined by
the function $Q_\alpha(t)$
and does not otherwise depend on the particular correlation function
$\rho_\alpha(t)$. In the following section, we discuss several examples in order to
illustrate the concept.

\section{Examples}\label{sec4}

In most cases, the asymptotic optimal designs for the regression model
(\ref{1.1})
have to be found numerically; explicit solutions are only possible for
very simple models.
In this section, we consider models
with one or two parameters.

\subsection{Optimal designs for linear models with one parameter}

Consider the linear regression model with $p=1$, that is,
$y(t)=\theta f(t)+\varepsilon(t)$ ($\theta\in\mathbb{R}$). In this case, the problem of minimizing the
asymptotic variance-covariance of the least-squares estimate reduces to the
minimization of the function
\begin{eqnarray*}
\Psi_\alpha (p)=\frac{\int f^2(t)Q_\alpha({\int p(u)
\,\mathrm{d}u}/{p(t)})p(t) \,\mathrm{d}t}{ (\int f^2(t)p(t)\, \mathrm{d}t )^2}\int p(t)
\,\mathrm{d}t.
\end{eqnarray*}
in the class of all non-negative functions $p(t)$ on the interval $[-T,T]$.
Note that we have represented the density $\phi$ by $p/\int p(x)\,\mathrm{d}x$,
which simplifies the calculation of the directional
derivatives in the following discussion.
Because $Q_\alpha(t)$ is strictly convex on $(0,\infty)$, it follows from
Theorem 3.1 in Bickel and Herzberg (\citeyear{BickelHerzberg1979}) that a minimizer,
say~$p^*(t)$, exists
and that $\phi^*(t)=p^*(t)/\int p^*(t)\,\mathrm{d}t$ is the asymptotic optimal
density. For the minimizing
function $p^*$, we obtain
\begin{eqnarray*}
\frac{\partial}{\partial\epsilon}\Psi_\alpha\bigl(p^*+\epsilon(p-p^*)\bigr)
\bigg|_{\epsilon=0}\ge0
\end{eqnarray*}
for all non-negative functions $p$ on the interval $[-T,T]$.
Consequently, the
asymptotic optimal density should satisfy $\int p^*(t)\,\mathrm{d}t=1$ and
%
\begin{eqnarray}\label{eq:opt_dens_eq}
\int\bigl(f^2(t)\bigl(H_\alpha\bigl(1/p^*(t)\bigr)-\mu\bigr) +\tilde\tau\bigr)\bigl(p(t)-p^*(t)\bigr)\,\mathrm{d}t\ge0
\end{eqnarray}
for all non-negative functions $p$ on the interval $[-T,T]$,
where the function \mbox{$H_\alpha{}\dvtx{} (0, \infty) \to\mathbb{R}^+$} is given by
%
\begin{eqnarray}\label{eq:mu}
H_\alpha(t)&=&Q_\alpha(t)-tQ'_\alpha(t)=
\cases{\displaystyle\frac{1+\alpha}{1-\alpha}/t^\alpha,&\quad$0<\alpha<1$,\cr
2/t,&\quad$\alpha=1$,
}
\nonumber
\\
\mu&=&2\frac{\int f^2(t)Q_\alpha(1/p^*(t))p^*(t)\,\mathrm{d}t}{ \int f^2(t)p^*(t)\,\mathrm{d}t },
\\
\tilde\tau&=&\int f^2(t)Q_\alpha(1/p^*(t))p^*(t)\, \mathrm{d}t +\int
f^2(t)Q'_\alpha(1/p^*(t)) \,\mathrm{d}t.\nonumber
\end{eqnarray}
First, assume that $0<\alpha<1$. Note that the
function $H_\alpha$ is strictly decreasing with $H_\alpha(+0)=\infty
$, $H_\alpha(\infty)=0$
and that its inverse is given by
\[
H^-_\alpha(t)=\biggl(\frac{1+\alpha}{t(1-\alpha)} \biggr)^{1/\alpha}.
\]
Hence the solution of
\eqref{eq:opt_dens_eq} has the form
%
\begin{eqnarray}\label{eq:opt_density}
\hspace*{-10pt}p^*(t)=
\cases{
\displaystyle\frac{1}{H^{-1}_\alpha\bigl(\mu-\tau/f^2(t)\bigr)}
 \cr
 \quad=\biggl( \displaystyle\frac{1-\alpha}{1 +\alpha}
\bigl(\mu-\tau/f^2(t) \bigr)\biggr)^{ 1 /\alpha},&\quad$\mu-\tau/f^2(t)\ge0$,\cr
0,&\quad\mbox{otherwise},
}
\end{eqnarray}
where $\mu$ is defined by \eqref{eq:mu} and
%
\begin{eqnarray}\label{eq:tau}
\tau=\int f^2(t)Q_\alpha\bigl(1/p^*(t)\bigr)p^*(t)\,\mathrm{d}t+\int f^2(t)Q'_\alpha
\bigl(1/p^*(t)\bigr)\,\mathrm{d}t.
\end{eqnarray}
Note that $\tau$ is a solution of the equation $\int p(t)\,\mathrm{d}t=1$.
Indeed, multiplying $f^2(t)H_\alpha(1/p^*(t))\equiv\mu f^2(t) -\tau$ by
$p^*(t)$ and integrating with respect to $t$ yields
\begin{eqnarray*}
\int f^2(t)H_\alpha\bigl(1/p^*(t)\bigr)p^*(t)\,\mathrm{d}t
=\int\bigl(\mu f^2(t) -\tau\bigr)p^*(t)\,\mathrm{d}t
=\mu\int f^2(t) p^*(t)\,\mathrm{d}t-\tau.
\end{eqnarray*}
Now, the definition of $H_\alpha(t)$ and $\mu$
gives
\begin{eqnarray*}
&&\int f^2(t)Q_\alpha\bigl(1/p^*(t)\bigr)p^*(t)\,\mathrm{d}t
-\int f^2(t)Q'_\alpha\bigl(1/p^*(t)\bigr)\,\mathrm{d}t
\\
&&\quad=2\int f^2(t)Q_\alpha\bigl(1/p^*(t)\bigr)p^*(t)\,\mathrm{d}t-\tau,
\end{eqnarray*}
which yields \eqref{eq:tau}.
Consequently, we have proven the following result.

\begin{theorem}
Assume that the correlation function $\rho_\alpha(t)$ is an element
of the Cauchy, Mittag-Leffler or SVF family. Then, for the one-parameter linear
regression model, the asymptotic optimal design
exists, is absolute continuous with respect to the Lebesgue measure and
has the density $p^*(t)$ defined in \eqref{eq:opt_density},
where $\mu$ and $\tau$ are given by \eqref{eq:mu} and~\eqref{eq:tau},
respectively.
\end{theorem}

\begin{figure}[t]

\includegraphics{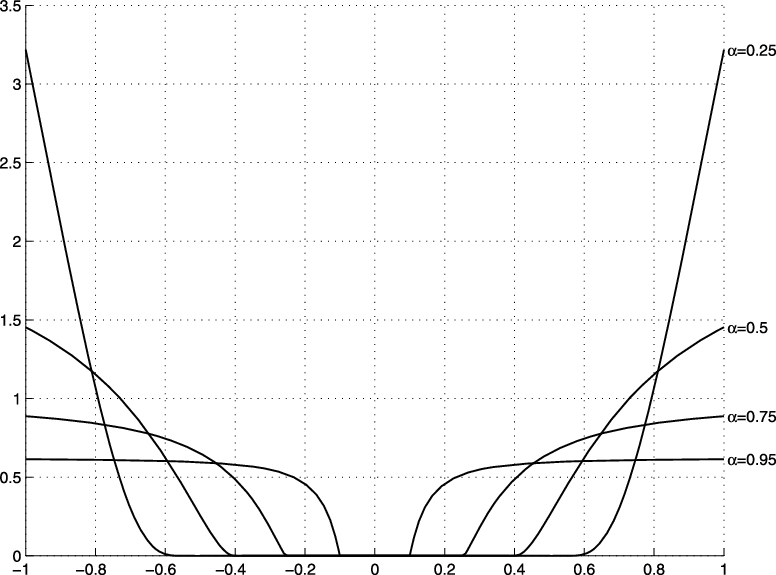}

\caption{Asymptotic optimal design densities for the linear
regression through the origin, $T=1$.} \label{fig1}
\end{figure}

We now consider two special cases, which are of particular importance. If
$p=1$ and $f(t)\equiv1$, then we obtain the location model and the
asymptotic optimal density is the uniform density, that is,
%
\begin{eqnarray}\label{eq:uni}
p^*(t)=
\cases{
\displaystyle\frac{1}{2T}, &\quad$|t|\le T$,\vspace*{2pt}\cr
0,&\quad\mbox{otherwise.}
}
\end{eqnarray}
Similarly, in the linear regression through the origin, we have
$p=1$, $f(t)\equiv t$ and the asymptotic optimal density is given
by
\begin{eqnarray*}
p(t)=
\cases{
0, &\quad$|t|\le\sqrt{\tau/\mu}$,\cr
\biggl(\displaystyle\frac{1-\alpha}{1+\alpha}\bigl(\mu-\tau/t^2\bigr)\biggr)^{1/\alpha},
&\quad$\sqrt{\tau/\mu}\le|t|\le T$,\vspace*{2pt}\cr
0, &\quad\mbox{otherwise},
}
\end{eqnarray*}
where
\begin{eqnarray*}
\mu=2\frac{\int t^2 p^{1+\alpha}(t)\,\mathrm{d}t}{ (1-\alpha)\int t^2p(t)\,\mathrm{d}t},\qquad
\tau=\int t^2 p^{1+\alpha}(t)\,\mathrm{d}t
\end{eqnarray*}
and $\alpha$ is the parameter of the correlation function.
The above formulas are given for $0<\alpha<1$. For $\alpha=1$ and
$f(t)= t$,
the asymptotic optimal density is the uniform density~\eqref{eq:uni}.
The optimal densities for the parameters $\alpha=1/4, 1/2, 3/4, 0.95$
and $T=1$ are displayed in Figure~\ref{fig1}.
The parameters $\mu$ and $\tau$ and the efficiency of uniform design
are shown
in Table~\ref{tab:tab1}. We observe that the uniform design is rather inefficient
for small values
of the parameter $\alpha$. The uniform design has
a reasonable efficiency only if $\alpha$ is close to~1.

%

%
\begin{table}[t]
\caption{Parameters of the asymptotic optimal design density for
the linear regression through the origin and
the efficiency of uniform design \protect\eqref{eq:uni}, $T=1$ (the optimal
density for $\alpha=1$ is $p(t)=1/2$, $-1\leq t \leq1$)}\label{tab:tab1}
\begin{tabular*}{\textwidth}{@{\extracolsep{\fill}}ld{2.2}ccc@{}}
\hline
$\alpha$
&\multicolumn{1}{c}{$\mu$}
&$\tau$
&$\sqrt{\tau/\mu}$
&$\operatorname{eff}_{\mathrm{uni}}$
\\
\hline
0.05& 2.34& 1.06& 0.67& 0.40\\
0.25& 3.19& 0.96& 0.55& 0.59\\
0.50& 4.32& 0.70& 0.40& 0.78\\
0.75& 6.84& 0.44& 0.25& 0.93\\
0.95& 24.78& 0.25& 0.10& 0.99\\
\hline
\end{tabular*}
\end{table}

It is worthwhile to mention that the asymptotic optimal designs derived
thus far depend sensitively on the parameter $\alpha$, which is usually not
available before the experiment. Because misspecification of this parameter can result in a substantial
loss of efficiency of the optimal design, we propose the construction of robust designs
which are less sensitive with respect to such misspecifications.
More precisely, we denote by $p^*_\alpha(t)$ the optimal density design for parameter $\alpha$.
Following Dette (\citeyear{Dette1995}) or
M\"{u}ller and P\'{a}zman (\citeyear{PazmanMuller1998}), a~robust
version of the optimality criterion is of the form
\begin{eqnarray*}
\Psi_\mathcal{A} (p)
=\min_{\alpha\in\mathcal{A}}\operatorname{eff}(p,\alpha)=
\min_{\alpha\in\mathcal{A}} \frac{\Psi_\alpha(p^*_\alpha)}{\Psi_\alpha(p)},
\end{eqnarray*}
where $p_\alpha^*$ is the optimal design for the correlation function
$\rho_\alpha$
and $\mathcal{A}$ is set of possible $\alpha$ values
specified by the experimenter. A design maximizing
$\Psi_\mathcal{A}$ is called \textit{standardized maximin optimal}.
Numerical optimization of this function for the set
$\mathcal{A}=\{0.1,0.2,\ldots,0.9\}$ shows that standardized maximin
optimal design has a density which can be approximated by the
function
%
\begin{eqnarray}
p^*_\mathcal{A} (t)=(5.7275t^2-1.16963-3.0264t^4)_+,
\label{eq:mm_des}
\end{eqnarray}
which is close to the optimal density $p^*_\alpha$ for $\alpha=0.44$.
In Table~\ref{tab:tab5}, we show the efficiency of this design for
various values of $\alpha$.
We observe that the design $p^*_\mathcal{A} $ is very efficient for all
elements in the set $\mathcal{A} $.

\subsection{Linear regression}

Consider the case $p=2$, $f_1(t)=1$, $f_2(t)=t$, which corresponds to
the linear regression
model. In this case, the asymptotic variance-covariance matrix is
proportional to
\begin{eqnarray*}
\Psi_\alpha(p)=
\pmatrix{
1 & \displaystyle\int t p(t)\,\mathrm{d}t\cr
\displaystyle\int t p(t)\,\mathrm{d}t
& \displaystyle\int t^2 p(t)\,\mathrm{d}t
}^{-1}
R(p)
\pmatrix{ & \displaystyle\int t p(t)\,\mathrm{d}t\cr
\displaystyle\int t p(t)\,\mathrm{d}t
& \displaystyle\int t^2 p(t)\,\mathrm{d}t
}
^{-1},
\end{eqnarray*}
where
\begin{eqnarray*}
R(p)=
\pmatrix{
\displaystyle\int Q_\alpha\bigl(1/p(t)\bigr)p(t)\,\mathrm{d}t
& \displaystyle\int t Q_\alpha\bigl(1/p(t)\bigr)p(t)\,\mathrm{d}t
\cr
\displaystyle\int t Q_\alpha\bigl(1/p(t)\bigr)p(t)\,\mathrm{d}t
& \displaystyle\int t^2Q_\alpha\bigl(1/p(t)\bigr)p(t)\,\mathrm{d}t
}
.
\end{eqnarray*}
For a symmetric density, this matrix is diagonal and
\begin{eqnarray*}
\Psi_\alpha(p)=\operatorname{diag}\biggl(\int Q_\alpha\bigl(1/p(t)\bigr)p(t)\,\mathrm{d}t,
\frac{\int t^2Q_\alpha(1/p(t))p(t)\,\mathrm{d}t}{(\int t^2 p(t)\,\mathrm{d}t)^2}\biggr).
\end{eqnarray*}
Consequently, the optimal symmetric design for estimating the slope in
the linear regression
has the density \eqref{eq:opt_density},
where $\mu$ and $\tau$ are defined in \eqref{eq:mu} and \eqref
{eq:tau} (this follows from the fact that
the element in position $(2,2)$ of the matrix
$\Psi_\alpha(p)$ corresponds to the optimality criterion for the
linear regression
through the origin). Numerical results indicate that the optimal
symmetric design for estimating the slope is optimal in the class of
all (not
necessarily symmetric) designs.

%

\begin{table}
\caption{Efficiency of the standardized maximin optimal design
$p^*_\mathcal{A} $ defined by \protect\eqref{eq:mm_des} in the linear regression through the origin
(the correlation structure is given by the SVF family with parameter
$\alpha$)}\label{tab:tab5}
\begin{tabular*}{\textwidth}{@{\extracolsep{\fill}}lccccccccc@{}}
\hline
$\alpha$
& 0.1
& 0.2
& 0.3
& 0.4
& 0.5
& 0.6
& 0.7
& 0.8
& 0.9
\\
\hline
$\operatorname{eff}(p^*,\alpha)$& 0.84& 0.92& 0.97& 0.99& 0.99& 0.97& 0.94&
0.89& 0.84\\
\hline
\end{tabular*}
\end{table}

\begin{figure}

\includegraphics{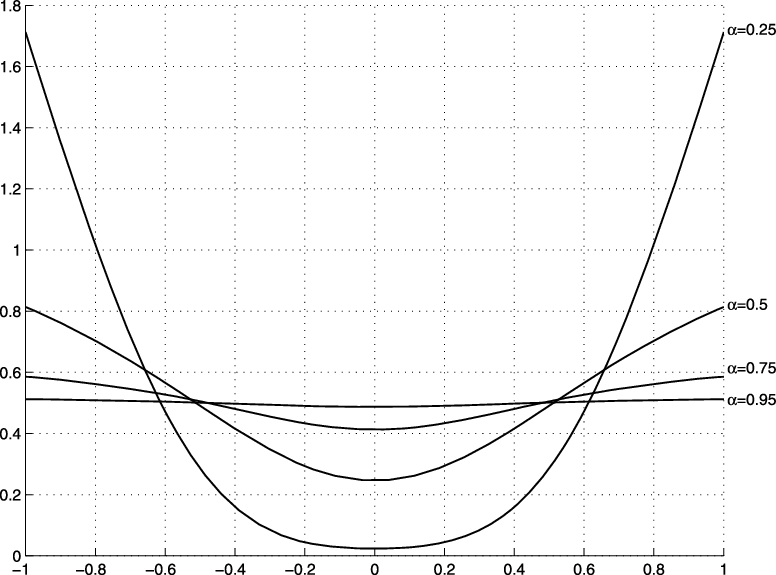}

\caption{Asymptotic D-optimal design densities for the linear regression,
$T=1$.}\label{figd}
\end{figure}

The D-optimal designs for the linear regression model have to be
determined numerically in all cases.
Some D-optimal design densities corresponding to the parameters $\alpha
=1/4,1/2,3/4,0.95$ and $T=1$ are displayed
in Figure~\ref{figd}.
%

\subsection{Comparison of optimal designs under long- and short-range
dependence}

It is of some interest to compare the asymptotic optimal designs under
short- and long-range dependence.
For this purpose, we again consider the linear regression model with no
intercept. Bickel and Herzberg (\citeyear{BickelHerzberg1979})
discussed the correlation function $\rho_\lambda(t)=\mathrm{e}^{-\lambda|t|}$. The asymptotic
optimal designs are given by
\begin{eqnarray*}
p(t)=
\cases{
0, &\quad$|t|\le\sqrt{\tau/\mu}$,\cr
\displaystyle\frac{1}{H^-(\mu-\tau/t^2)},
&\quad$\sqrt{\tau/\mu}\le|t|\le T$,\cr
0, &\quad\mbox{otherwise},
}
\end{eqnarray*}
%
\begin{figure}

\includegraphics{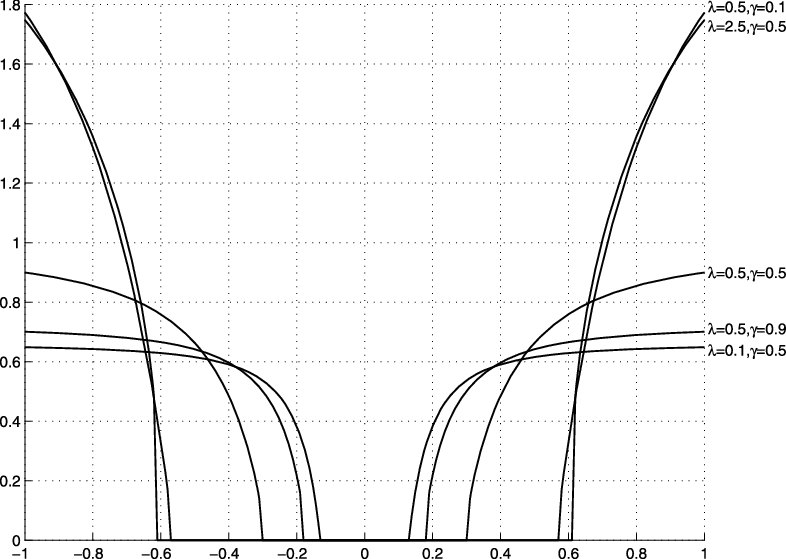}%
\vspace*{-3pt}
\caption{Asymptotic optimal design densities for the
linear regression through the origin, where the correlation function is
given by $\rho_\lambda(t)=\mathrm{e}^{-\lambda|t|}$.} \label{fig2}\vspace*{-6pt}
\end{figure}
where the quantities $\mu, \tau$ are defined by
\begin{eqnarray*}
\mu&=&\frac{1}{2\gamma}+2\frac{\int f^2(t)Q_\alpha(1/p^*(t))p^*(t)\,\mathrm{d}t}{ \int f^2(t)p^*(t)\,\mathrm{d}t },
\\
\tau&=&\frac{1}{2\gamma}\int f^2(t)p^*(t)\,\mathrm{d}t +\int f^2(t)Q_\alpha
\bigl(1/p^*(t)\bigr)p^*(t)\,\mathrm{d}t+
\int f^2(t)Q_\alpha'\bigl(1/p^*(t)\bigr)\,\mathrm{d}t,
\end{eqnarray*}
respectively (see Bickel, Herzberg and Schilling (\citeyear{BickelHerzbergSchilling1981}))
and depend on the parameters $\lambda$ and $\gamma$ defined in \eqref
{eq:corr_fun_r}.
Some of these designs are shown in Figure~\ref{fig2}, while the
relevant parameters
%
\begin{table}
\caption{Parameters of the asymptotic optimal design density for
the linear regression through the origin, where the correlation
function is
given by $\rho_\lambda(t)=\mathrm{e}^{-\lambda|t|}$ (the last column of the
table shows the efficiency of the uniform design \protect\eqref{eq:uni})}\label{tab:tab1a}
\vspace*{-3pt}
\begin{tabular*}{\textwidth}{@{\extracolsep{\fill}}lccccc@{}}
\hline
$\lambda$
&$\gamma$
&$\mu$
&$\tau$
&$\sqrt{\tau/\mu}$
&$\operatorname{eff}_{\mathrm{uni}}$
\\
\hline
0.5& 0.5& 3.41& 0.32& 0.30& 0.89\\
0.5& 0.1& 9.82& 3.23& 0.57& 0.63\\
0.5& 0.9& 2.38& 0.08& 0.18& 0.97\\
0.1& 0.5& 12.70& 0.22& 0.13& 0.99\\
2.5& 0.5& 1.45& 0.54& 0.61& 0.57\\
\hline
\end{tabular*}
\end{table}
are given in Table~\ref{tab:tab1a}, which also contains the efficiency
of the uniform design. We observe
that -- in contrast to the case of long-range dependence -- the uniform
design is rather efficient, provided either that the parameter
$\lambda$ is not too large or that $\gamma$ is not too small.

We now compare asymptotic optimal designs derived under the assumption
of a long-range dependence
with asymptotic optimal designs under short-range dependence. In Table
\ref{tab:tab3}, we show the efficiency
of a design derived under the assumption of short-range dependence, in
the situation where the ``true''
correlation structure is a member of the SVF family. We observe that
the loss of efficiency
is only substantial if the parameter $\alpha$ is small. The opposite
situation is displayed in Table
\ref{tab:tab4}, which shows the efficiency of the
asymptotic optimal design under long-range dependence (from the
SVF family), but the ``true'' correlation structure is in fact of
exponential type. Again, the asymptotic optimal designs
derived under the long-range dependence are rather efficient, except
when the parameter $\alpha$ is very small.

\begin{table}
\caption{Efficiency of the asymptotic
optimal design density for the correlation function
$\rho_\lambda(t)=\mathrm{e}^{-\lambda|t|}$ in
the linear regression through the origin, while the ``true''
correlation function belongs
to the SVF family}\label{tab:tab3}
\begin{tabular*}{\textwidth}{@{\extracolsep{\fill}}lcccccc@{}}
\hline
& $\alpha$
& 0.05
& 0.25
& 0.50
& 0.75
& 0.95\\
\hline
$\lambda=0.5$& $\gamma=0.5$& 0.62& 0.82& 0.96& 1.00& 0.97\\
$\lambda=0.5$& $\gamma=0.1$& 0.81& 0.97& 0.99& 0.89& 0.77\\
$\lambda=0.5$& $\gamma=0.9$& 0.53& 0.73& 0.90& 0.99& 1.00\\
$\lambda=0.1$& $\gamma=0.5$& 0.50& 0.70& 0.88& 0.98& 1.00\\
$\lambda=2.5$& $\gamma=0.5$& 0.81& 0.97& 0.98& 0.89& 0.77\\
\hline
\end{tabular*}\vspace*{-8pt}
\end{table}

\begin{table}
\caption{Efficiency of asymptotic optimal design density for a
long-range dependence error structure in the linear regression through
the origin, while the ``true'' correlation function is given by
$\rho_\lambda(t)=\mathrm{e}^{-\lambda|t|}$}\label{tab:tab4}
\begin{tabular*}{\textwidth}{@{\extracolsep{\fill}}lccccc@{}}
\hline
$\lambda$
& 0.5
& 0.5
& 0.5
& 0.1
& 2.5
\\
$\gamma$
& 0.5
& 0.1
& 0.9
& 0.5
& 0.5\\
\hline
$\alpha=0.05$& 0.19& 0.40& 0.15& 0.15& 0.35\\
$\alpha=0.25$& 0.69& 0.94& 0.59& 0.58& 0.93\\
$\alpha=0.50$& 0.94& 0.98& 0.87& 0.86& 0.98\\
$\alpha=0.75$& 1.00& 0.88& 0.98& 0.98& 0.85\\
$\alpha=0.95$& 0.95& 0.73& 0.99& 1.00& 0.68\\
\hline
\end{tabular*}\vspace*{-8pt}
\end{table}

\section*{Acknowledgements}
The support of the Deutsche
Forschungsgemeinschaft (SFB 475,
``Komplexit\"{a}tsreduktion in multivariaten
Datenstrukturen'' SFB 823 ``Statistik nichtlinearer dynamischer Prozesse'')
is gratefully acknowledged. The work of H.~Dette was also supported in
part by an NIH Grant award IR01GM072876:01A1 and by the BMBF project SKAVOE.
The second author, N. Leonenko, was supported in part by the EPSRC Grant RCMT 199, the Australian
Research Grant DP 0559807 and the Welsh Institute of
Mathematics and Computational Sciences. The authors would also like to
thank two referees for their constructive comments on an earlier version
of this paper.

\printhistory

\end{document}